\documentclass{article}

\RequirePackage[OT1]{fontenc}
\RequirePackage[
noinfoline,amsthm,amsmath
]{imsart}

\usepackage{graphicx}
\usepackage{latexsym,amsmath}
\usepackage{amsmath,amsthm,amscd}
\usepackage{amsfonts}
\usepackage[psamsfonts]{amssymb}
\usepackage{enumerate}

\usepackage{multirow}

\usepackage{url}

\usepackage{tocvsec2}

\usepackage[final]
{showkeys}

\usepackage{epstopdf}
\usepackage{float}
\usepackage[small,hang]{caption}

\usepackage{hyperref}

\startlocaldefs
\numberwithin{equation}{section}

\theoremstyle{plain} 
\newtheorem{theorem}{Theorem}[section]

\theoremstyle{definition} 

\theoremstyle{definition} 

\newtheorem*{ex*}{Example}
\theoremstyle{remark} 

\theoremstyle{remark} 
\newtheorem{remark}[theorem]{Remark}
\newtheorem*{remark*}{Remark}
\theoremstyle{definition} 

\numberwithin{equation}{section}

\renewcommand{\bar}[1]{{\overline{#1}}}

\newcommand{\sign}{\operatorname{sign}}

\newcommand{\tr}{\operatorname{tr}}

\renewcommand{\th}{\operatorname{th}}

\renewcommand{\th}{\theta}

\newcommand{\al}{\alpha}

\newcommand{\si}{\sigma}
\newcommand{\Si}{\Sigma}
\newcommand{\ka}{\kappa}

\newcommand{\be}{\beta}
\newcommand{\De}{\Delta}

\newcommand{\La}{\Lambda}
\newcommand{\vpi}{\varphi}

\renewcommand{\P}{\operatorname{\mathsf{P}}} 
\newcommand{\E}{\operatorname{\mathsf{E}}}
\newcommand{\Var}{\operatorname{\mathsf{Var}}}
\newcommand{\Cov}{\operatorname{\mathsf{Cov}}}

\newcommand{\R}{\mathbb{R}}

\newcommand{\vp}{\varepsilon}

\renewcommand{\le}{\leqslant}
\renewcommand{\ge}{\geqslant}

\endlocaldefs

\begin{document}

\begin{frontmatter}

\title{An asymptotically optimal transform of Pearson's 
correlation statistic
}
\runtitle{Optimal 
transform of Pearson's  
statistic}

%

\begin{aug}
\author{\fnms{Iosif} \snm{Pinelis}
\ead[label=e1]{ipinelis@mtu.edu}}
\runauthor{Iosif Pinelis}


\address{Department of Mathematical Sciences\\
Michigan Technological University\\
Houghton, Michigan 49931, USA\\
E-mail: \printead[ipinelis@mtu.edu]{e1}}
\end{aug}

\begin{abstract}
It is shown that for any correlation-parametrized model of dependence and any given significance level $\al\in(0,1)$, 
there is an asymptotically optimal transform of Pearson's correlation statistic $R$, for which the generally leading error term for the normal approximation vanishes 
for all values $\rho\in(-1,1)$ of the correlation coefficient. 
This general result is then applied to the bivariate normal (BVN) model of dependence and to what is referred to in this paper as the SquareV model.  
In the BVN model, Pearson's $R$ turns out to be asymptotically optimal for a rather unusual significance level $\al\approx0.240$, whereas Fisher's transform $R_F$ of $R$ is asymptotically optimal for the limit significance level $\al=0$. 
In the SquareV model, Pearson's $R$ is asymptotically optimal for a still rather high significance level $\al\approx0.159$, whereas Fisher's transform $R_F$ of $R$ is not asymptotically optimal for any $\al\in[0,1]$.
Moreover, it is shown that in both the BVN model and the SquareV model, 
the transform optimal for a given value of $\al$ is in fact asymptotically better than $R$ and $R_F$ in wide ranges of values of the significance level, including $\al$ itself. 
Extensive computer simulations for the BVN and SquareV models of dependence are presented, which suggest that, for sample sizes $n\ge100$ and significance levels $\al\in\{0.01,0.05\}$, the mentioned asymptotically optimal transform of $R$ generally outperforms both Pearson's $R$ and Fisher's transform $R_F$ of $R$, the latter appearing generally much inferior to both $R$ and the asymptotically optimal transform of $R$ in the SquareV model. 
\end{abstract}

  
%


%


\end{frontmatter}

\settocdepth{chapter}

\tableofcontents 

\settocdepth{subsubsection}

\theoremstyle{plain} 
\numberwithin{equation}{section}



\section{
Introduction}\label{intro} 

A statistic closely related to Pearson's $R$ is commonly known as the Fisher $z$ transform, defined by the formula 
\begin{equation}
	R_F:=\tanh^{-1}(R)=\frac12\ln\frac{1+R}{1-R}.
\end{equation}
An advantage of using $R_F$ (as opposed to $R$) in making statistical inferences about the true correlation coefficient $\rho$ is usually ascribed to its variance-
stabilizing property in normal populations; see e.g.\ Fisher \cite{fish15}, Gayen \cite{gayen51}, and Hotelling \cite{hot53}; that is, $n\Var R_F\to1$ for all $\rho\in(-1,1)$ as $n\to\infty$ \big(as opposed to $n\Var R \to(1-\rho^2)^2$\big) whenever 
the underlying distribution is bivariate normal. 
Everywhere here, $n$ denotes the sample size. 

In his discussion of Hotelling's paper \cite{hot53}, Kendall
provides heuristics suggesting that such variance stabilization of the distribution of a statistic may often result in it being closer to normality. 
Namely, if an approximate constancy of the variance of a statistic were the same as an approximate constancy of its distribution itself, and if the distribution is close to normality at least for one value of the parameter (say, $\rho$, as in the present case), then it would be close to normality for all values of $\rho$. 
For normal populations and large enough sample sizes, the Fisher $z$ transform indeed brings the distribution of the correlation statistic closer to normality, 
and it is especially effective for values of $\rho$ far from $0$. 
However, it is well known 
(see e.g.\ \cite{BG78,nonlinear-publ}) 
that the closeness of the distribution of a statistic to normality is usually mainly determined, not by the variance, but by the third moments of the appropriately standardized statistic. 


In this paper, we shall see that for a general and most common class of models of dependence, including the bivariate normal (BVN) model, and
for each given significance level $\al\in(0,1)$ there is a certain transform $\Psi_\al(R)$ of Pearson's statistic $R$ that assures the \emph{vanishing} of the generally leading term of the asymptotics of the probability that an approximately standardized version of the statistic $\Psi_\al(R)$ exceeds the standard normal critical value 
\begin{equation}\label{eq:z_al}
	z_\al:=\Phi^{-1}(1-\al); 
\end{equation}
here, as usual, $\Phi$ is the standard normal cumulative distribution function (cdf) and $\Phi^{-1}$ is its inverse; unless otherwise specified, all the asymptotics here are for large sample sizes $n$.  
Thus, the 
transform $\Psi_\al(R)$ of $R$ is asymptotically optimal: its distribution is asymptotically the closest to normality exactly at the critical value. 

Once this optimality result is obtained for the general class of models of dependence, the rest of the paper is devoted to detailed analysis of the optimal transform $\Psi_\al(R)$ of $R$ in the BVN model and another specific model of dependence, referred to in this paper as the SquareV model. 

We shall 
see that, in the BVN model, the mentioned family $(\Psi_\al(R))_{\al\in(0,1)}$ of transforms of Pearson's statistic $R$ includes $R$ itself: namely, 
\begin{equation}\label{eq:R=Psi(R)}
	R=\Psi_{\al_P}(R), 
\end{equation}
where 
\begin{equation*}
	\al_P:=1-\Phi(1/\sqrt2)\approx0.240. 
\end{equation*}
Thus, Pearson's statistic $R$ is asymptotically optimal for a significance level $\al$ of about $24\%$, but such a significance level is rather unusual in statistical practice. 

As for Fisher's transform $R_F$ of $R$, we shall see that, again in the BVN model,  
\begin{equation}\label{eq:R_F}
	R_F=\Psi_0(R):=\lim_{\al\downarrow0}\Psi_\al(R), 
\end{equation}
which means that Fisher's transform is asymptotically optimal for the significance level 
\begin{equation*}
	\al_F:=0. 
\end{equation*}
Now one might explain the fact that for the usually rather small significance levels, such as $0.05$ or $0.01$, Fisher's statistic $R_F$ is asymptotically closer to normality than Pearson's statistic $R$ by noting that the significance level $\al_F=0$ (for which $R_F$ is asymptotically optimal) is closer to $0.05$ and especially to $0.01$ than the significance level $\al_R\approx0.240$ (for which $R$ is asymptotically optimal). 

As for the SquareV model, there Pearson's $R$ is asymptotically optimal for a still rather high significance level $\al\approx0.159$, whereas Fisher's transform $R_F$ of $R$ is not asymptotically optimal for any $\al\in[0,1]$.

It should be noted that in both the BVN model and the SquareV model, 
the transform optimal for a given value of $\al$ is in fact asymptotically better than $R$ and $R_F$ in wide ranges of values of the significance level, including $\al$ itself. 

We have also conducted extensive computer simulations for the BVN and SquareV models of dependence, which suggest that, for sample sizes $n\ge100$ and significance levels $\al\in\{0.01,0.05\}$, the mentioned asymptotically optimal transform of $R$ generally outperforms both Pearson's $R$ and Fisher's transform $R_F$ of $R$, the latter appearing generally much inferior to both $R$ and the asymptotically optimal transform of $R$ in the SquareV model.

The rest of the paper is organized as follows. 

In Section~\ref{chib} we present an asymptotic expansion for statistics that are general smooth nonlinear functions of the sample mean of iid random vectors in $\R^p$. This expansion, which may be viewed as a far-reaching refinement of the delta method, is a special case of results by Bhattacharya and Ghosh \cite{BG78}. For Berry--Esseen-type bounds for general nonlinear statistics, see e.g. \cite{chen07,nonlinear-publ}.  

In Section~\ref{R and z asymp}, the mentioned asymptotic expansion is specialized for the cases of Pearson's correlation statistic $R$ and its smooth enough transforms. A key observation there is that the main term of the asymptotic for such a transform of $R$ differs from the corresponding main term for $R$ itself only by a comparatively simple expression involving the first two derivatives of the transform function $\psi$. This allows one to obtain, for any correlation-parametrized model of dependence and for any given significance level $\al\in(0,1)$, a rather simple second-order ordinary differential equation (ODE) for the optimal transform function $\psi$ that makes the main term of the asymptotic for the asymptotically optimal transform of $R$ vanish for all values $\rho\in(-1,1)$ of the correlation coefficient. This ODE can be explicitly solved for a number of models of dependence, including the important BVN model and models with a linear dependence of the correlation parameter. The mentioned SquareV model is a model with such a linear dependence. 

The BVN model is considered in detail in Section~\ref{BVN}, where the corresponding results of the mentioned computer simulations are also presented and discussed. 

A similar treatment of the SquareV model is given in Section~\ref{4points}. 

Section~\ref{concl} is a summary of the results of this paper.

\section{Asymptotic expansions for smooth nonlinear statistics}\label{chib} 

Let 
\begin{equation*}
	V,V_1,V_2,\dots
\end{equation*}
be independent identically distributed (iid) zero-mean random vectors in $\R^p$ with 
$\E\|V\|^3<\infty$, where $p$ is a natural number and $\|\cdot\|$ is the Euclidean norm in $\R^p$, which latter will be identified, as usual, with the space of all $p\times1$ column matrices. 
Assume also that the Cram\'er-type condition $\limsup_{\|t\|\to\infty}|\E\exp(it^T V)|
<1$ is satisfied, where $i$ is the imaginary unit and ${}^T$ denotes the transposition, in this case of a column matrix $t\in\R^p$; for this Cram\'er-type condition to hold, it is enough  
that, for some natural $k$, the $k$-fold convolution of the distribution of $V$ have a nonzero absolutely continuous component. 
Let $\Si$ stand for the covariance matrix of $V$: 
\begin{equation}\label{eq:Si}
	\Si:=\E VV^T. 
\end{equation}

Let $f\colon\R^p\to\R$ be a function which is twice continuously differentiable in a neighborhood of $0\in\R^p$ and such that $f(0)=0$. Let $L$ and $H$ denote, respectively, the gradient vector and the Hessian matrix of the function $f$ at $0$, so that 
\begin{equation}\label{eq:L,H}
	f'(0)(v)=L^T v\quad\text{and}\quad f''(0)(v,v)=v^T Hv
\end{equation}
for all $v\in\R^p$. 
Since $V$ is assumed to be zero-mean, one has $\E L^T V=0$. 
Introduce now 
\begin{equation}\label{eq:si}
	\si:=\sqrt{\E(L^T V)^2}, 
\end{equation}
which will be assumed to be nonzero, so that 
\begin{equation*}
	\La:=\frac{L^T V}\si
\end{equation*}
is a well-defined r.v., with zero mean and unit variance. 
Consider the r.v.\   
\begin{equation*}
	T_n:=\frac{\sqrt n}\si\,f(\bar V),  
\end{equation*}
where of course $\bar V:=\frac1n\,\sum_1^n V_i$. 
Then, by Theorem~2 of the paper \cite{BG78} by Bhattacharya and Ghosh, 
\begin{equation}\label{eq:BG}
	\sup_{z\in\R}|\P(T_n\le z)-\Psi_{3,n}(z)|=o\Big(\frac1{\sqrt n}\Big), 
\end{equation}
where 
%
%
\begin{align}
	\Psi_{3,n}(z)&:=\Phi(z)+\frac{\De(z)}{\sqrt n}, \\ 
		\De(z)&:=-\Big[\Big(\frac{\E\La^3}6+a_3\Big)\,(z^2-1)+a_1\Big]\,\vpi(z) \label{eq:De sum} \\ 
	&=(Az^2+B)\,\vpi(z),   \label{eq:A,B}  
\end{align}
$\Phi$ and $\vpi$ denote, as usual, the distribution and density functions of $N(0,1)$, and $a_1$, $a_2$, $A$, and $B$ are constants depending only on $L$, $H$, $\Si$, $\si$, and $\E\La^3$ (but not on $z$ or $n$): 
\begin{align}
	a_1&:=\frac1{2\si}\,\tr H\Si, \label{eq:a_1} \\ 
		a_3&:=\frac1{4\si^3}\,(L^T\Si L-\si^2)\tr(H\Si)+\frac1{2\si^3}\,L^T\Si H\Si L, \label{eq:a_3} \\ 
	A&:=-\Big(\frac{\E\La^3}6+a_3\Big), \\
	B&:=-A-a_1,  
\end{align}
with $\tr$ denoting the trace of a matrix. 

\begin{remark}\label{rem:O(sqrt n)}
If the condition $\E\|V\|^3<\infty$ is strengthened to $\E\|V\|^4<\infty$, then $o\big(\frac1{\sqrt n}\big)$ in \eqref{eq:BG} can be replaced by $O\big(\frac1n\big)$. 
\end{remark}

\section{Asymptotics for the Pearson statistic and its transforms} \label{R and z asymp}
Let $(Y,Z),(Y_1,Z_1),\dots,(Y_n,Z_n)$ be independent identically distributed random points in $\R^2$ with a correlation coefficient $\rho\in(-1,1)$ and $\E(Y^6+Z^6)<\infty$. Pearson's sample correlation coefficient based on the observations $(Y_1,Z_1),\dots,(Y_n,Z_n)$ is defined by the formula 
\begin{equation}\label{eq:R}
	R:=R_n:=\frac{\bar{YZ}-\bar Y\,\bar Z}{\sqrt{\bar{Y^2}-\bar Y^2}\sqrt{\bar{Z^2}-\bar Z^2}},   
\end{equation}
where 
$\bar Y:=\frac1n\sum_1^n Y_i$, $\bar Z:=\frac1n\sum_1^n Z_i$, $\bar{Y^2}:=\frac1n\sum_1^n Y_i^2$, $\bar{Z^2}:=\frac1n\sum_1^n Z_i^2$, and $\bar{YZ}:=\frac1n\sum_1^n Y_iZ_i$; 
%
let 
$R$ take an arbitrarily assigned value in the interval $[-1,1]$ if the denominator of the ratio in \eqref{eq:R} is $0$. 

Let us assume that $Y$ and $Z$ are each standardized, that is, zero-mean and unit-variance. This assumption does not diminish generality, because $R$ is invariant with respect to affine transformations $Y_i\mapsto a+b\,Y_i$ and $Z_i\mapsto c+d\,Z_i$ of the $Y_i$'s and $Z_i$'s, for any real constants $a,b,c,d$ such that $b>0$ and $d>0$. 

Observe that 
\begin{equation}\label{eq:=f}
	R-\rho=f(\bar V), 
\end{equation} 
where 
\begin{equation}\label{eq:V_R}
	V:=(Y,Z,Y^2-1,Z^2-1,YZ-\rho) 
\end{equation}
and 
\begin{equation}\label{eq:f_R}
	f(v):=f_\rho(v):=\frac{\rho+v_5-v_1v_2}{\sqrt{1+v_3-v_1^2}\sqrt{1+v_4-v_4^2}}-\rho  
\end{equation}
if $v=(v_1,\dots,v_5)\in\R^5$ is such that $1+v_3-v_1^2>0$ and $1+v_4-v_4^2>0$; 
otherwise, let $f(v):=
0$. 
In this case, $L=(0,0,-\frac\rho2,-\frac\rho2,1)$, whence 
\begin{equation}\label{eq:si,La}
\si=\sqrt{\E\big(YZ-\tfrac\rho2\,(Y^2+Z^2)\big)^2}\quad\text{and}\quad 
	\La=\frac{YZ-\tfrac\rho2\,(Y^2+Z^2)}\si. 
\end{equation}

As noted in \cite{nonlinear-publ}, the condition $\si=0$ is equivalent to the following exceptional situation: there exists some $\ka\in\R$ such that the random point $(Y,Z)$ lies almost surely on the union of the two straight lines through the origin with slopes $\ka$ and $1/\ka$ (for $\ka=0$, these two lines should be understood as the two coordinate axes in the plane $\R^2$). 

It will be assumed in what follows that the random point $(Y,Z)$ is such that $\si$ is never $0$. 
 
Then it is easy to check that all the conditions on $f$ and $V$ stated in Section~\ref{chib} are satisfied, with $p=5$; 
in particular, $f(0)=0$.

Let now $\De_R(z)$ denote $\De(z)$ defined by \eqref{eq:De sum} with $f$ as in \eqref{eq:f_R}. 

Further, letting 
\begin{equation}\label{eq:mu_ij}
	\mu_{ij}:=\E Y^iZ^j, 
\end{equation}
one has 
\begin{equation}\label{eq:si_R}
   \si=\frac12\,\sqrt{\rho^2 \left(\mu _{04}+2 \mu _{22}+\mu _{40}\right)-4 \rho \left(\mu _{13}+\mu
   _{31}\right)+4 \mu _{22}},   	
\end{equation}
\begin{multline}\label{eq:ELa3}
	\si^3\E\La^3=-\frac{\rho^3}{8} \left(\mu _{06}+3 \mu _{24}+3 \mu _{42}+\mu_{60}\right)
	+6 \rho^2 \left(\mu _{15}+2 \mu _{33}+\mu _{51}\right) \\ 
   -12 \rho\left(\mu _{24}+\mu _{42}\right)+8 \mu _{33}, 
\end{multline}
and 
\begin{multline}\label{eq:De_R}
	\frac{96\si^3}{\vpi(z)} \De_R(z) = \tilde\De_R(z)\\ 
:=16 \big[(z^2-1) (6 \mu_{12} \mu_{21}-\mu_{33})+3 \si^2 z^2 (\mu_{13}+\mu_{31})\big]  \hfill \\ 
-12 \rho \big[(z^2-1) (4 \mu_{03} \mu_{21}+4
   \mu_{12} \mu_{30}+8 \mu_{12}^2-2 \mu_{13} \mu_{31}+\mu_{13}^2+8 \mu_{21}^2-2 \mu_{24}+\mu_{31}^2-2 						\mu_{42}) \\ 
\hfill    +\si^2 \big((2z^2+1) (\mu_{04}+\mu_{40})+(4 z^2-2) \mu_{22}
   \big)\big] \\ 
+12 \rho^2 (z^2-1) \big[2 \mu_{03} (3 \mu_{12}+\mu_{30})+\mu_{04} (\mu_{13}-\mu_{31})+10 \mu_{12} 								\mu_{21}-\mu_{13} \mu_{40}-\mu_{15} \\ 
\hfill    +6 \mu_{21} \mu_{30}+\mu_{31} \mu_{40}-2 \mu_{33}-\mu_{51}\big] \\ 	
-\rho^3 (z^2-1) \big[24 \mu _{03} \mu _{21}+12 \mu _{03}^2-6 
   \mu_{04} \mu_{40}+3 \mu_{04}^2-2 \mu_{06} \hfill \\
\hfill    +24 \mu_{12} \mu_{30}+12 \mu_{12}^2+12 \mu_{21}^2-6 \mu_{24}+12 \mu_{30}^2+3 \mu_{40}^2-6 \mu_{42}-2 			\mu_{60}\big]. 
\end{multline}

More generally, let now 
\begin{equation*}
\psi\colon(-1,1)\to\R 	
\end{equation*}
be a twice continuously differentiable function whose derivative $\psi'$ does not vanish at any point of the interval $(-1,1)$. Let then 
\begin{equation}\label{eq:g=}
	g(v):=g_\rho(v):=\frac{\psi\big(f(v)+\rho\big)-\psi(\rho)}{\psi'(\rho)}
\end{equation}
for $v\in\R^5$, with $f$ as defined in \eqref{eq:f_R}. 
Then, in view of \eqref{eq:=f}, 
\begin{equation}\label{eq:=g}
\frac{\psi(R)-\psi(\rho)}{\psi'(\rho)}=g(\bar V), 
\end{equation}
Note that 
$g(0)=f(0)=0$ and $g'(0)=f'(0)$, so that $\si$ and $\La$ for the function $g$ are the same as in \eqref{eq:si,La} (given there for the function $f$). 
Thus, all the conditions stated in Section~\ref{chib} are satisfied with $g$ in place of $f$. 

Let now $\De_{\psi(R)}(z)$ denote $\De(z)$ defined by \eqref{eq:De sum} with $g$ as in \eqref{eq:g=} in place of $f$. 
Then \eqref{eq:BG} will hold with 
\begin{equation}\label{eq:tau_n}
	\tau_{\psi,n}:=\frac{\psi(R)-\psi(\rho)}{\psi'(\rho)\si/\sqrt n}
\end{equation}
in place of $T_n$ and $\De_{\psi(R)}(z)$ in place of $\De(z)$.  
One may note here that $\tau_{\psi,n}$ may be considered an asymptotically standardized version of $\psi(R)$.

A key observation is that 
\begin{equation}\label{eq:differ}
	\De_{\psi(R)}(z)=\De_R(z)-\frac{\psi''(\rho)}{2\psi'(\rho)}\,\si z^2\,\vpi(z). 
\end{equation}

To begin using this observation, let us refer to any family $(P_\rho)_{\rho\in(-1,1)}$ of distributions of the random pair $(Y,Z)$ in $\R^2$ parametrized by the correlation coefficient $\rho$ of $(Y,Z)$ as a \emph{correlation-parametrized model (CP) of dependence}. 
The CP condition seems quite natural for parametric models of dependence. Indeed, let a real parameter $\th$ represent the strength of the dependence between $Y$ and $Z$. Then one should usually expect the correlation coefficient $\rho$ to be a strictly increasing continuous function $g$ of $\th$: $\rho=g(\th)$. Replacing then $\th$ by $g^{-1}(\rho)$, one obtains a re-parametrization with $\rho$ as the new parameter. 

In this regard, one may recall the formula 
\begin{equation}\label{eq:cov=}
	\Cov(Y,Z)=\iint_{\R^2}a_{Y,Z}(y,z)\,dy\,dz
\end{equation}
for the covariance of $Y$ and $Z$, where $a_{Y,Z}$ is the association function of r.v.'s $Y$ and $Z$, given by the formula 
\begin{equation*}
	a_{Y,Z}(y,z):=\P(Y>y,Z>z)-\P(Y>y)\P(Z>z)
\end{equation*}
for all $(y,z)\in\R^2$. 
Since $\P(Y>y)=\P(-Y<-y)$, formula \eqref{eq:cov=} can be rewritten as 
\begin{equation}\label{eq:cov=,F}
	\Cov(Y,Z)=\iint_{\R^2}[F_{Y,Z}(y,z)-F_Y(y)\,F_Z(z)]\,dy\,dz, 
\end{equation}
where $F_{Y,Z}$ is the joint cdf of the random pair $(Y,Z)$, and $F_Y$ and $F_Z$ are the corresponding marginal cdf's. 

Suppose now that we have a dependence model $(\P_\th)$, where $\th$ is a strength of the association/dependence parameter, so that the association function $a_{\th;Y,Z}$ of the pair $(Y,Z)$ with respect to the probability measure $\P_\th$ is increasing in $\th$ on the average in the sense that 
the integral in \eqref{eq:cov=} with $a_{\th;Y,Z}$ in place of $a_{Y,Z}$ is increasing in $\th$. Suppose also that the $Y$- and $Z$-marginals of the distribution of the pair $(Y,Z)$ with respect to $\P_\th$ do not depend on $\th$. Then the correlation coefficient $\rho$ of $(Y,Z)$ will be an increasing function of $\th$. 
 
A more specific, but still rather general way to construct a CP model is as follows. By Sklar's theorem (see e.g.\ \cite[Theorem 2.3.3.]{nelsen06}), 
\begin{equation*}
	F_{Y,Z}(y,z)=C(F_Y(y),F_Z(z))
\end{equation*}
for some copula $C$ and all $(y,z)\in\R^2$; recall that a copula can be defined as the joint cdf of a random pair with values in the unit square $[0,1]^2$ whose marginals are uniform on the interval $[0,1]$. 
Let now $(C_\th)$ be any family of copulas increasing in $\th$; a large number of such families can be found in \cite{nelsen06}. Fix the marginal cdf's $F_Y$ and $F_Z$, and for each value of the parameter $\th$ let $F_{\th;Y,Z}(y,z):=C_\th(F_Y(y),F_Z(z))$, again for all $(y,z)\in\R^2$. Then the correlation coefficient $\rho$ corresponding to the joint cdf $F_{\th;Y,Z}$ will be an increasing function of $\th$. 

In view of \eqref{eq:mu_ij}, in any correlation-parametrized model of dependence and for any given real $z\ne0$,  
%
the expressions in \eqref{eq:si_R} and \eqref{eq:De_R} for  
$\si$ and $\tilde\De_R(z)$ 
will depend on $\rho$ only. 
Then, by the key observation \eqref{eq:differ}, the condition $\De_{\psi(R)}(z)=0$ can be rewritten as the second-order ordinary differential equation (ODE)
\begin{equation}\label{eq:ODE}
	\frac{\psi''(\rho)}{\psi'(\rho)}=h_z(\rho)  
\end{equation}
for the function $\psi$, where 
\begin{equation}\label{eq:h_z}
	h_z(\rho):=\frac{\tilde\De_R(z)}{48\si^4 z^2}.   
\end{equation}
Solving now ODE \eqref{eq:ODE} with the natural initial conditions
\begin{equation}\label{eq:init}
	\psi(0)=0\quad\text{and}\quad\psi'(0)=1, 
\end{equation}
we have 
\begin{equation}\label{eq:psi'=}
	\psi'(\rho)=\exp\int_0^\rho dr\, h_z(r)
\end{equation}
and 
\begin{equation}\label{eq:psi=}
		\psi(\rho)=\psi_z(\rho):=\int_0^\rho dr\,\exp\int_0^r ds\, h_z(s)
\end{equation}
for $\rho\in(-1,1)$; in \eqref{eq:psi'=} and \eqref{eq:psi=}, 
we use the common convention 
$\int_0^s:=-\int_s^0$ for $s<0$. 
Thus, we obtain 

\begin{theorem}\label{th:general}
In any correlation-parametrized model of dependence and for any given real $z\ne0$, 
the generally leading error term for the normal approximation for $\psi_z(R)$ vanishes:  
\begin{equation}\label{eq:=0,z}
	\De_{\psi_z(R)}(z)=0 
\end{equation}
for all $\rho\in(-1,1)$. 
\end{theorem}

Letting now 
\begin{equation}\label{eq:Psi:=}
	\Psi_\al:=\psi_{z_\al},
\end{equation}
we can rewrite \eqref{eq:=0,z} as 
\begin{equation*}
	\De_{\Psi_\al(R)}(z_\al)=0 
\end{equation*}
for all $\al\in(0,1)$, with $z_\al=\Phi^{-1}(1-\al)$, as defined in \eqref{eq:z_al}. 

One may note here that for a rather large class of models of dependence the functions $h_z$ will be rational, and hence, according to \eqref{eq:psi'}, $\psi'_z$ will be an elementary, closed-form function. This class of models with rational functions $h_z$ includes 
the bivariate normal model and models with linear dependence of the joint cdf $F_{\th;Y,Z}$ on $\th$. In particular, the class of models with linear dependence of $F_{\th;Y,Z}$ on $\th$ contains Farlie's model \cite{far60}. 

We shall consider the bivariate normal model and a particular simple model with linear dependence of $F_{\th;Y,Z}$ on $\th$ in the following sections, to compare the performance of Pearson's $R$ itself, its Fisher transform $R_F$, and the asymptotically optimal transform $\Psi_\al(R)$ of $R$ in non-asymptotic settings, for specific sample sizes. 


\section{Bivariate normal model (BVN)}\label{BVN}

\subsection{Asymptotically optimal transform \texorpdfstring{$\Psi_\al(R)$}{} in the BVN model}\label{Psi-BVN}

Here it is assumed that the random point $(Y,Z)$ has the bivariate normal (BVN) distribution with zero means, unit variances, and an arbitrary correlation coefficient $\rho\in(-1,1)$. Then the expressions for $\tilde\De_R$ and $\si$, and thus for $h_z(\rho)$, in formulas \eqref{eq:De_R}, \eqref{eq:si_R}, and \eqref{eq:h_z} can be greatly simplified. 

Indeed, in this case the pair $(Y,Z)$ equals $(Y,\rho\, Y+\sqrt{1-\rho^2}\,Y_1)$ in distribution, whence, by \eqref{eq:mu_ij},   
\begin{equation*}
	\mu_{ij}=\sum _{k=0}^j \binom{j}{k} \rho^{k} \left(1-\rho^2\right)^{(j-k)/2} m(i+k)\,m(j-k)
\end{equation*}
for all $i,j=0,1,\dots$ and 
\begin{equation*}
	m_j:=\E Y^j, 
\end{equation*}
so that $(m_0,\dots,m_6)=(1,0,1, 0, 3, 0, 15)$. 
As the result, ODE \eqref{eq:ODE} becomes  
\begin{equation}\label{eq:ODE N}
	\frac{\psi''(\rho)}{\psi'(\rho)}=
	p_z\,\frac{-2\rho}{1-\rho ^2},  
\end{equation}
where 
\begin{equation}\label{eq:p_z}
	p_z:=\frac1{2z^2}-1. 
\end{equation}
ODE \eqref{eq:ODE N} is easily solved, yielding 
\begin{equation}\label{eq:psi'}
	\psi'(\rho)=(1 - \rho^2)^{p_z}
\end{equation}
and 
\begin{equation}\label{eq:psi N}
	\psi(\rho)=\psi_z(\rho):=\int_0^\rho(1 - r^2)^{p_z}\,dr
	=\rho  \  _2F_1\big(\tfrac{1}{2},-p_z;\tfrac{3}{2};\rho ^2\big), 
\end{equation}
where $_2F_1$ is the ordinary hypergeometric function, given by the formula 
\begin{equation*}
	_2F_1(a,b;c;x)=\sum_{k=0}^\infty\frac{(a)_k(b)_k}{(c)_k}\,\frac{x^k}{k!}
\end{equation*}
for $x$ with $|x|<1$, where $(q)_k:=\prod_{i=0}^{k-1}(q+i)$ is the Pochhammer symbol. 
The last equality in \eqref{eq:psi N} can be obtained by expanding the integrand $(1 - r^2)^{p_z}$ into the Maclaurin series in powers of $r$ and then integrating the series term-wise. 

Thus, recalling \eqref{eq:Psi:=}, we see that in the bivariate normal case the transform $\Psi_{\al}(R)=\psi_{z_\al}(R)$ of Pearson's $R$ with $\psi_z$ as in \eqref{eq:psi N} is asymptotically optimal for any given significance level $\al\in(0,1)$. 

In particular, choosing
\begin{equation*}
z=1/\sqrt2\approx0.707,	
\end{equation*}
we have $p_z=0$. Hence, in view of the integral expression in \eqref{eq:psi N}, $\psi_z(\rho)=\rho$, so that we have   
\eqref{eq:R=Psi(R)}, confirming that the family $(\Psi_\al(R))_{\al\in(0,1)}$ of transforms of Pearson's statistic $R$ includes $R$ itself. 

On the other hand, letting $z\to\infty$, we have $p_z\to-1$, so that, using again the integral expression in \eqref{eq:psi N} (and, say, the dominated convergence theorem), we see 
\begin{equation}\label{eq:psi_infty}
	\psi_z(\rho)\underset{z\to\infty}\longrightarrow\psi_\infty(\rho):=\frac12\ln\frac{1+\rho}{1-\rho} 
\end{equation}
for all $\rho\in(-1,1)$, thus confirming \eqref{eq:R_F}. 

In view of formula \eqref{eq:psi N}, the calculation of values of the functions $\psi_z$ or, equivalently, of the functions $\Psi_\al$ mainly reduces to the calculation of values of the hypergeometric function $_2F_1$. In general, this hypergeometric function is not elementary. However, there are a number of highly efficient ways to compute values of $_2F_1$. It takes only about $1.7\times10^{-5}$ sec 
on an average to compute a value of $\psi_2(\rho)$ (on a standard computer), which may be compared with the corresponding execution time of about $0.45\times10^{-5}$ sec for Fisher's $\psi_\infty(\rho)=\frac12\ln\frac{1+\rho}{1-\rho}$. 
Therefore and because usually in statistical practice the value of the transform $\Psi_\al(R)=\psi_{z_\al}(R)$ of the statistic $R$ needs to be computed only once, the use of the hypergeometric function $_2F_1$ should not cause any complications. 

Also, according to \eqref{eq:psi'}, the derivarive $\psi'$ of the function $\psi=\psi_z$ is a simple elementary expression, which makes it easy to obtain various analytical properties of $\psi_z$. For instance, using the special l'Hospital-type rule for monotonicity (see e.g.\ \cite[Proposition 4.1]{pin06}), we can immediately see that the ratio $\psi_{z_1}(\rho)/\psi_{z_2}(\rho)$ is decreasing in $\rho^2$ for any real $z_1$ and $z_2$ such that $0<|z_1|<|z_2|$. In particular, it follows that the ratio of  $\psi_z(\rho)$ to Fisher's $\psi_\infty(\rho)=\frac12\ln\frac{1+\rho}{1-\rho}$ is decreasing in $\rho^2$ for any real $z\ne0$. One may also note that the values  
\begin{equation*}
	\psi_z(\pm1)=\pm\frac{\sqrt{\pi }\, \Gamma (p_z+1)}{2\, \Gamma \left(p_z+3/2\right)}
\end{equation*}
at the endpoints of the interval $[-1,1]$ 
are finite for all real $z\ne0$, in contrast with Fisher's limit values $\psi_\infty(\pm(1-))=\pm\infty$; here $p_z$ is as defined in \eqref{eq:p_z}. 
It is also clear that $\psi_z(\rho)$ is odd in $\rho$, for each $z\ne0$. 
 

\subsection{The transform \texorpdfstring{$\Psi_\al(R)$}{} in the BVN model is asymptotically better than \texorpdfstring{$R$}{} and \texorpdfstring{$R_F$}{} in wide ranges of values of the significance level, including \texorpdfstring{$\al$}{} itself} 
\label{other-al-BVN}

According to Theorem~\ref{th:general}, for any given correlation-parametrized model of dependence and any given significance level $\al\in(0,1)$, the transform $\Psi_\al(R)$ of $R$ is asymptotically optimal for all $\rho\in(-1,1)$ as the sample size $n$ goes to $\infty$. In particular, for any given significance level $\al\in(0,1)$, the transform $\Psi_\al(R)$ of $R$ is asymptotically better than both $R$ itself and its Fisher transform $R_F$. In fact, $\Psi_\al(R)$ is asymptotically better than $R$ and $R_F$ for rather wide ranges of values (say $\be$) of the significance level; of course, these ranges include the value $\al$ itself. 
Indeed, one can see that in the BVN model 
\begin{equation*}
\De_{\psi(R)}(z)=\frac\rho2\,(2 z^2-1)-\frac{z^2}{2} \frac{(1-\rho ^2)\psi''(\rho )}{\psi'(\rho)}. 
\end{equation*}
In particular, 
\begin{equation*}
\De_R(z)=\frac\rho2\,(2 z^2-1), 
\end{equation*} 
\begin{equation*}
\De_{R_F}(z)=\De_{\psi_\infty(R)}(z)=-\frac\rho2, 
\end{equation*}
where $\psi_\infty$ is as defined in \eqref{eq:psi_infty}, and 
\begin{equation*}
\De_{\Psi_\al(R)}(z)=\De_{\psi_{z_\al}(R)}(z)=\frac\rho2\Big(\frac{z^2}{z_\al^2}-1\Big). 
\end{equation*} 
So, 
$|\De_{\Psi_\al(R)}(z_\be)|<|\De_{R_F}(z_\be)|$ for $\rho\ne0$ if $0<z_\be<z_\al\sqrt2$. That is,
the transform $\Psi_\al(R)$, which is asymptotically optimal for the given significance level $\al\in(0,1)$, will still be asymptotically better than Fisher's transform $R_F$ for any significance level $\be\in(0,1)$ such that $0<z_\be<z_\al\sqrt2$. For instance, if $\al=0.05$, then $\Psi_\al(R)$ will be asymptotically better than $R_F$, not just for the significance level $\al=0.05$, but for any significance level $\be\in(0.01000,0.5)$ -- because $0<z_\be<z_{0.05}\sqrt2$ for all $\be\in(0.01000,0.5)$. Similarly, if $\al=0.01$, then $\Psi_\al(R)$ will be asymptotically better than $R_F$ for any significance level $\be\in(0.00050,0.5)$. 

As for the comparison of the asymptotically optimal transform $\Psi_\al(R)$ with $R$ itself, we can similarly see that, for instance, if $\al=0.05$, then $\Psi_\al(R)$ will be asymptotically better than $R$, not just for the significance level $\al=0.05$, but for any significance level $\be\in(0,0.17912)$; if $\al=0.01$, then $\Psi_\al(R)$ will be asymptotically better than $R$ for any significance level $\be\in(0,0.16933)$. 

\subsection{Comparison of asymptotically optimal transform \texorpdfstring{$\Psi_\al(R)$}{} in the BVN model with \texorpdfstring{$R$}{} and \texorpdfstring{$R_F$}{} for finite sample sizes} 
\label{finite-BVN}

According again to Theorem~\ref{th:general}, for any given correlation-parametrized model of dependence and any given significance level $\al\in(0,1)$, the transform $\Psi_\al(R)$ of $R$ is asymptotically optimal for all $\rho\in(-1,1)$ as the sample size $n$ goes to $\infty$. However, one may ask: how does this asymptotically optimal transform of $R$ compare in performance with $R$ itself and the Fisher transform of $R$ for finite sample sizes $n$? 

To address this question for the BVN model, we have done the following. For each value $\al\in\{0.01,0.05\}$ of the significance level, each value $\rho\in\{0,0.1,0.5,0.9\}$ of the correlation coefficient, and each value $n\in\{10,10^2,10^3,10^4\}$ of the sample size, $N:=10^6$ samples of size $n$ from the BVN distribution with zero means, unit variances, and correlation coefficient $\rho\in(-1,1)$ were simulated using the command \texttt{RandomVariate[BinormalDistrib\-ution[], \{\}]} of the computing system 
Mathematica, on each of the 12 working in parallel processing cores of a computer workstation. Thus, altogether $12\times10^6\times2\times4\times(10+10^2+10^3+10^4)\approx10^{12}$ realizations of the BVN random pair $(Y,Z)$ were simulated, which took about $3.6$ hours. This suggests that Mathematica's simulation of the BVN distribution is highly effective. \big(In view of the symmetry $\Cov(Y,-Z)=-\Cov(Y,Z)$, negative values of $\rho$ have not been considered here.\big)

For each quadruple $$(\al,\rho,n,k)\in\{0.01,0.05\}\times\{0,0.1,0.5,0.9\}\times\{10,10^2,10^3,10^4\}\times\{1,\dots,12\}$$ 
(where $k$ indexes the 12 processing cores) and for each of the three transforms of $R$ \big($R$ itself, its Fisher transform $R_F$, and the asymptotically optimal transform $\Psi_\al(R)$ of $R$\big), the relative frequencies -- say $\hat\al_{\psi(R)}(\al,\rho,n,k)$ -- of the values of 
$\tau_{\psi,n}$ exceeding $z_\al$ were computed, where $\psi(\rho)\equiv\psi_{1/\sqrt2}(\rho)\equiv\rho$ for $R$, $\psi=\psi_\infty$ for $R_F$, $\psi=\psi_{z_\al}$ for $\Psi_\al(R)$; and $\tau_{\psi,n}$ is as defined in \eqref{eq:tau_n}. So, each of these relative frequencies $\hat\al_{\psi(R)}(\al,\rho,n,k)$ is an estimate of the significance level $\al$, with the corresponding estimates 
\begin{equation*}
\hat\vp_{\psi(R)}(\al,\rho,n,k):=\frac{\hat\al_{\psi(R)}(\al,\rho,n,k)}\al-1  	
\end{equation*}
of the relative errors 
of the approximation of $\al$ by $\P_\rho(\tau_{\psi,n}>z_\al)$, for each $\psi\in\{\psi_{1/\sqrt2},\psi_\infty,
\psi_{z_\al}\}$. 
Then, for each triple $(\al,\rho,n)\in\{0.01,0.05\}\times\{0,0.1,0.5,0.9\}\times\{10,10^2,10^3,10^4\}$ and each $\psi\in\{\psi_{1/\sqrt2},\psi_\infty,\psi_{z_\al}\}$, the mean 
\begin{equation*}
	\hat\vp_{\psi(R)}(\al,\rho,n):=\frac1{12}\sum_{k=1}^{12}\hat\vp_{\psi(R)}(\al,\rho,n,k)
\end{equation*}
and standard deviation 
\begin{equation*}
s_{\psi(R)}(\al,\rho,n):=\sqrt{
\frac1{12-1}
\sum_{k=1}^{12}\big(\hat\vp_{\psi(R)}(\al,\rho,n,k)-\hat\vp_{\psi(R)}(\al,\rho,n)\big)^2, 
}	
\end{equation*}
of the $12$ values 
\begin{equation}\label{eq:hateps}
	\hat\vp_{\psi(R)}(\al,\rho,n,1),\dots,\hat\vp_{\psi(R)}(\al,\rho,n,12)
\end{equation}
were computed. Thus, again for each triple $(\al,\rho,n)$ and each $\psi\in\{\psi_{1/\sqrt2},\psi_\infty,\psi_{z_\al}\}$, $\hat\vp_{\psi(R)}(\al,\rho,n)$ is the best available in this setting estimate of the relative error of the approximation of $\al$ by 
$\P_\rho(\tau_{\psi,n}>z_\al)$, and a reasonable estimate of the standard error of this estimate of the relative error of the approximation is 
\begin{equation*}
	\tilde s_{\psi(R)}(\al,\rho,n):=s_{\psi(R)}(\al,\rho,n)/\sqrt{12}. 
\end{equation*}


The values of $\hat\vp_{\psi(R)}(\al,\rho,n)$ and $\tilde s_{\psi(R)}(\al,\rho,n)$ are presented in Table~\ref{tab:tabBVN}, along with the corresponding values of $\tilde s_{\psi(R)}(\al,\rho,n)$, all of those values round to the precision of $3$ decimal digits. 

For instance, the entry $-0.236\pm0.00253$ in the first row of Table~\ref{tab:tabBVN}, under the heading $\hat\vp_R(\al,\rho,n)\pm\tilde s_R(\al,\rho,n)$, means that for the triple $(\al,\rho,n)=(0.01 , 0 , 10)$ and $\psi(R)=\psi_{1/\sqrt2}(R)=R$, we have $\hat\vp_{\psi(R)}(\al,\rho,n)=\hat\vp_R(\al,\rho,n)\approx-0.236$ and 
$\tilde s_{\psi(R)}=\tilde s_R(\al,\rho,n)\approx0.00253$. For the same triple $(\al,\rho,n)=(0.01 , 0 , 10)$, we have $\hat\vp_{R_F}(\al,\rho,n)=\hat\vp_{\psi_\infty(R)}(\al,\rho,n)\approx1.63$ and 
$\tilde s_{R_F}(\al,\rho,n)=\tilde s_{\psi_\infty(R)}(\al,\rho,n)\approx0.00463$, and also 
$\tilde s_{\Psi_\al(R)}(\al,\rho,n)=\tilde s_{\psi_{z_\al}(R)}(\al,\rho,n)\approx0.00439$. 

\begin{figure}[h]
	\centering
		\includegraphics[width=1.00\textwidth]{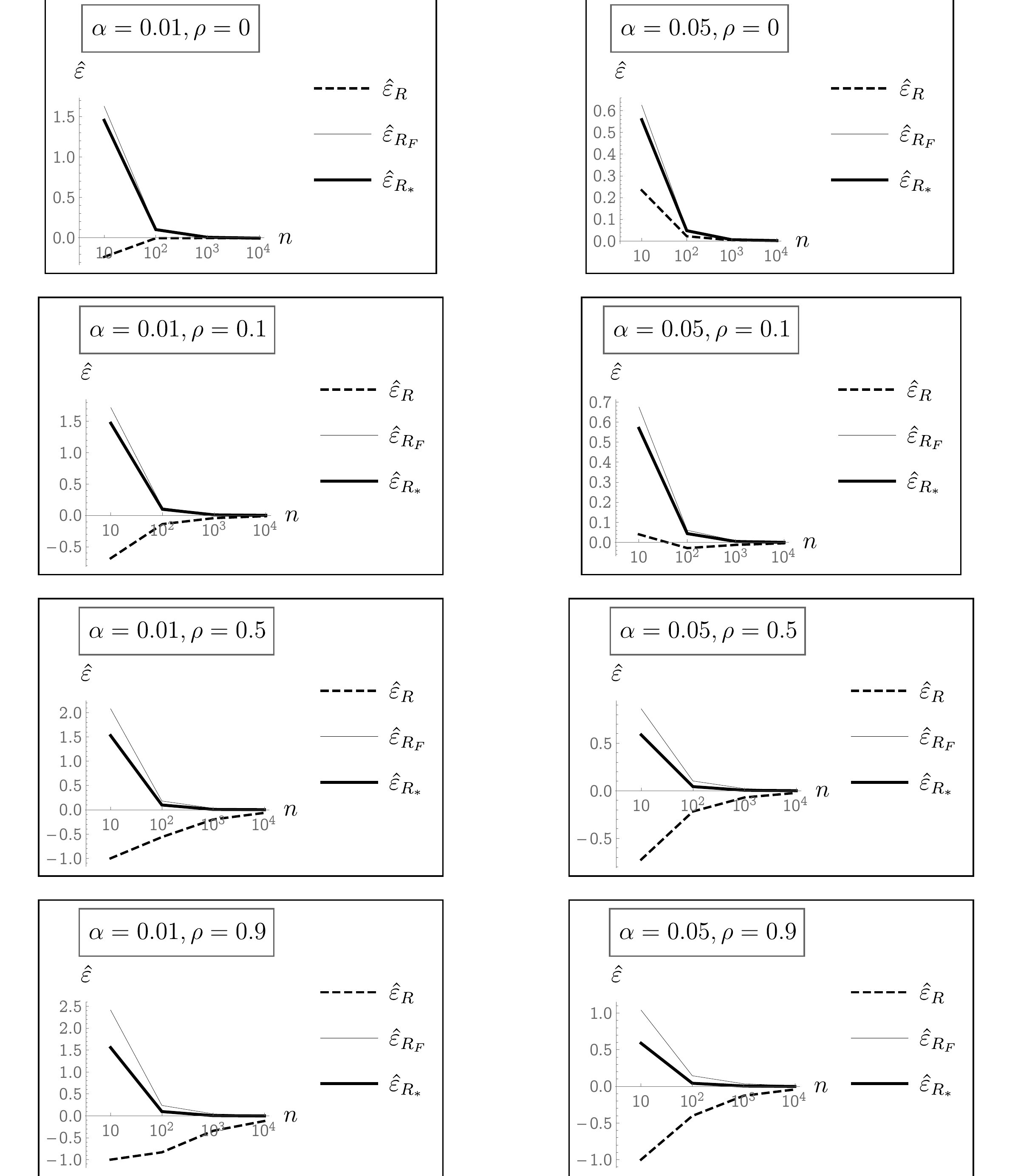}
	\caption{
	Means 
	of the 
estimates \eqref{eq:hateps} of the relative errors of the approximation of the significance level $\al$ for Pearson's $R$ (thin), Fisher's transform $R_F$ of $R$ (dashed), and the asymptotically optimal transform $R_*:=\Psi_\al(R)$ (thick) in the BVN model}	
	\label{fig:pic1BVN}
\end{figure}

For an easier grasp, graphs of the values of the means $\hat\vp_{\psi(R)}(\al,\rho,n)$ are given in Fig.~\ref{fig:pic1BVN}. The 8 pictures in Fig.~\ref{fig:pic1BVN} correspond to the pairs $(\al,\rho)\in\{0.01,0.05\}\times\{0,0.1,0.5,0.9\}$ and are labeled accordingly. Each of these 8 pictures contains, for the given pair $(\al,\rho)$ and for the three transforms -- $R$, $R_F$, and $R_*:=\Psi_\al(R)$ of $R$, the successively connected graphs of $\hat\vp_{\psi(R)}(\al,\rho,n)$ as functions of $n\in\{10,10^2,10^3,10^4\}$. 

\begin{figure}[h]
	\centering
		\includegraphics[width=1.00\textwidth]{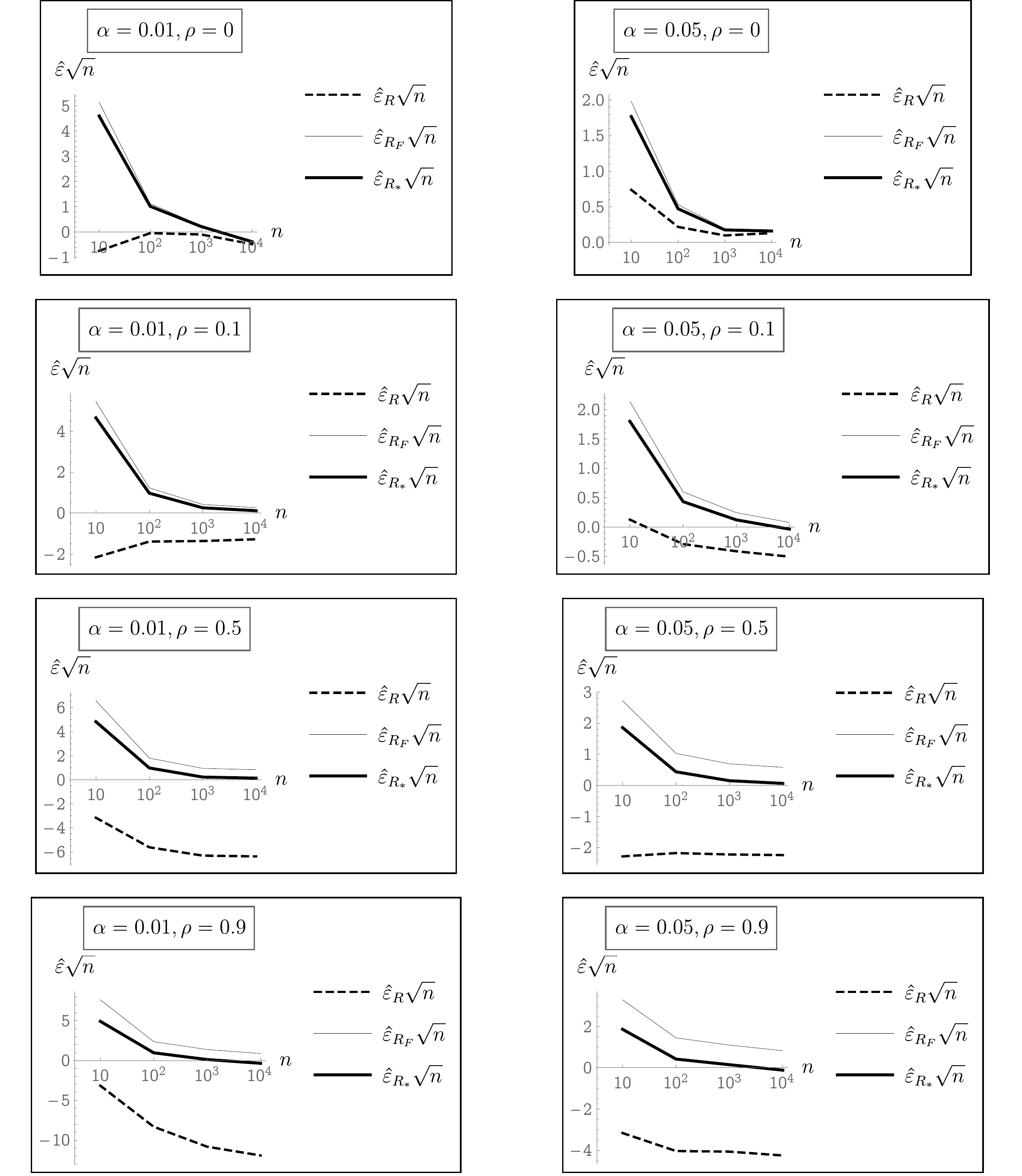}
		\caption{
	Means as in Fig.\ \ref{fig:pic1BVN} multiplied by $\sqrt n$ 
	} 
	\label{fig:pic2BVN}
\end{figure}

By Remark~\ref{rem:O(sqrt n)}, it should be expected that the means $\hat\vp_{\psi(R)}(\al,\rho,n)$ decrease in $n$ roughly as $1/\sqrt n$. Therefore, to see the differences between the three graphs in each of the 8 pictures more clearly for the larger values of $n$, the corresponding graphs of the values of  $\hat\vp_{\psi(R)}(\al,\rho,n)\,\sqrt n$ are given in Fig.~\ref{fig:pic2BVN}. 

Looking at Fig.~\ref{fig:pic2BVN} and back at Table~\ref{tab:tabBVN}, we see that for all considered triples $(\al,\rho,n)$ the asymptotically optimal transform $\Psi_\al(R)$ provides for a better normal approximation than Fisher's transform $R_F$ does. The advantage of $\Psi_\al(R)$ over $R_F$ is smaller for the smaller values of $\rho$ and $\al$, which should be understandable. 
Indeed, since all the three versions of the function $\psi$ (namely, $\psi_{1/\sqrt2}$, $\psi_\infty$, and $\psi_{z_\al}$) satisfy the same initial conditions \eqref{eq:init} at $\rho=0$, all the three corresponding transforms of $R$ ($R$ itself, $R_F$, and $\Psi_\al(R)$) are close to one another if $\rho$ is close to $0$. Also, the smaller the value of $\al$ is, the closer is the critical value $z_\al$ to Fisher's value $z_{0+}=\infty$. 

For the sample size $n=10$, the asymptotically optimal transform $\Psi_\al(R)$ appears to have little or no advantage over $R$ itself, which is also understandable, since the sample size $n=10$ is small, whereas the transform $\Psi_\al(R)$ is guaranteed to be optimal only for large enough sample sizes; actually, for $n=10$ all the three transforms perform poorly or, mostly, very poorly -- except for the rather strange case of $(\al,\rho,n)=(0.05,0.1,10)$, when $R$ performs surprisingly well. Also, the advantage of $\Psi_\al(R)$ over $R$ appears to be smaller for smaller values of $\rho$, a possible reason for which is that, as explained in the previous paragraph, all the three transforms of $R$ are close to one another if $\rho$ is close to $0$; in fact, $R$ appears to be somewhat better than $\Psi_\al(R)$ for $\rho=0$ (except when $\al=0.01$, $\rho=0$, and $n=10^4$) and for $\al=0.05$, $\rho=0.1$, and $n=10^2$. However, for $\rho\in\{0.5,0.9\}$ and $n\in\{10^2,10^3,10^4\}$, the transform $\Psi_\al(R)$ appears to be much better than $R$ itself. 

For the same values of $\rho$ and $n$ as in the last sentence, Fisher's transform $R_F$ also appears to be much better than $R$, but even in these cases, the advantage of $R_F$ over $R$ is significantly less than that of $\Psi_\al(R)$ over $R$. 

\section{SquareV model}\label{4points}

Here we shall consider the dependence model, which is the family $(P_\rho)_{-1<\rho<1}$ of distributions of the random pair $(Y,Z)$ on the vertices of the square $[-1,1]\times[-1,1]$ given by the following formulas: 
\begin{equation}\label{eq:SM1}
\begin{aligned}
	&\P_\rho\big((Y,Z)=(1,1)\big)=\P_\rho\big((Y,Z)=(-1,-1)\big)=\tfrac{1+\rho}4, \\ 
	&\P_\rho\big((Y,Z)=(1,-1)\big)=\P_\rho\big((Y,Z)=(-1,1)\big)=\tfrac{1-\rho}4.  
\end{aligned}	
\end{equation}
In other words, the distribution $P_{Y,Z}=P_{\rho;Y,Z}$ of $(Y,Z)$ under $\P_\rho$ is the mixture 
\begin{equation}\label{eq:SM2}
 \left\{
\begin{alignedat}{2}
&(1-\rho)P_{\vp_1,\vp_2}+\rho P_{\vp_1,\vp_1}&&\text{ if }\rho\ge0, \\ 
&(1+\rho)P_{\vp_1,\vp_2}-\rho P_{\vp_1,-\vp_1}&&\text{ if }\rho<0,  
\end{alignedat}
 \right. 
\end{equation}
where $\vp_1,\vp_2$ are independent Rademacher r.v.'s, with $\P(\vp_j=\pm1)=1/2$ for $j=1,2$. 
Then $\Cov(Y,Z)=\rho$ under $\P_\rho$, so that the use of the symbol $\rho$ to denote the parameter is consistent. 

The just described model of dependence will be referred to as the SquareV model, where ``V'' stands for ``vertices''. 

In view of the previously mentioned symmetry $\Cov(Y,-Z)=-\Cov(Y,Z)$, negative values of $\rho$ will not be further considered in this section. 

\subsection{Asymptotically optimal transform \texorpdfstring{$\Psi_\al(R)$}{} in the SquareV model}\label{Psi-4}

Using \eqref{eq:SM1} and \eqref{eq:SM2}, we obtain the following expressions for the joint moments of $(Y,Z)$ as defined in \eqref{eq:mu_ij}:   
\begin{align}
	\mu_{ij}&=\frac{1+\rho}{4}
   \left(1+(-1)^{i+j}\right) +\frac{1-\rho}{4} \left((-1)^i+(-1)^j\right)  \\
   &=(1-\rho)\,\frac{1+(-1)^i}2\,\frac{1+(-1)^j}2\,
   +\rho\,
   \frac{1+(-1)^{i+j}}2 
\end{align}
for all $i,j=0,1,\dots$ and all $\rho\in[0,1)$.

As the result, ODE \eqref{eq:ODE} becomes  
\begin{equation}\label{eq:ODE 4}
	\frac{\psi''(\rho)}{\psi'(\rho)}=
	q_z\,\frac{-2\rho}{1-\rho ^2},  
\end{equation}
where 
\begin{equation}\label{eq:q_z}
	q_z:=\frac1{3z^2}-\frac13.  
\end{equation}
We see that formulas \eqref{eq:ODE 4}--\eqref{eq:q_z} are rather similar to 
\eqref{eq:ODE N}--\eqref{eq:p_z}. 
Hence, quite similarly to \eqref{eq:psi'} and \eqref{eq:psi N}, 
here we have 
\begin{equation}\label{eq:psi'4}
	\psi'(\rho)=(1 - \rho^2)^{q_z}
\end{equation}
and 
\begin{equation}\label{eq:psi 4}
	\psi(\rho)=\psi_{4;z}(\rho):=\int_0^\rho(1 - r^2)^{q_z}\,dr
	=\rho  \  _2F_1\big(\tfrac{1}{2},-q_z;\tfrac{3}{2};\rho ^2\big). 
\end{equation}
Here the subscript $4$ in $\psi_{4;z}$ refers to the four points of the distribution of the random point $(Y,Z)$ in the SquareV model, currently under consideration; thus, one can distinguish between the function $\psi_{4;z}$ in \eqref{eq:psi 4} and the function $\psi_z$ in \eqref{eq:psi N}. 

Accordingly, recalling again \eqref{eq:Psi:=}, we see that in the SquareV model the transform 
\begin{equation*}
\Psi_{4;\al}(R)=\psi_{4;z_\al}(R)  	
\end{equation*}
of Pearson's $R$ 
is asymptotically optimal for any given significance level $\al\in(0,1)$. 

In particular, choosing
\begin{equation*}
z=1,	
\end{equation*}
we have $q_z=0$. Hence, in view of the integral expression in \eqref{eq:psi 4}, $\psi_{4;z}(\rho)=\rho$, so that the family $(\Psi_{4;\al}(R))_{\al\in(0,1)}$ of transforms of Pearson's statistic $R$ includes $R$ itself. 
More specifically, 
\begin{equation*}
	R=\Psi_{4;\al}(R)\quad\text{for}\quad \al=1-\Phi(1)\approx0.159. 
\end{equation*}

However, in order for Fisher's transform $R_F$ of $R$ to belong to the family $(\Psi_{4;\al}(R))_{\al\in(0,1)}$ of asymptotically optimal transforms of $R$ in the SquareV model, one would have to have $q_z=-1$ for some real $z$, which is impossible, because, in view of \eqref{eq:q_z}, $q_z$ is always greater than $-\frac13$. 
So, in contrast with the BVN model \big(where, according to \eqref{eq:R_F}, $R_F$ is asymptotically optimal in the limit case with $\al=0$ and $z_\al=\infty$\big), in the SquareV model Fisher's transform $R_F$ is not asymptotically optimal for any significance level $\al\in[0,1]$, even if the endpoints $\al=0$ and $\al=1$ are included as limit cases.  

\subsection{The transform \texorpdfstring{$\Psi_\al(R)$}{} in the SquareV model is asymptotically better than \texorpdfstring{$R$}{} and \texorpdfstring{$R_F$}{} in wide ranges of values of the significance level, including \texorpdfstring{$\al$}{} itself} 
\label{other-al-4}

This subsection is similar to Subsection~\ref{other-al-BVN}.
One can see that in the SquareV model 
\begin{equation*}
\De_{\psi(R)}(z)=\frac\rho{3\sqrt{1-\rho ^2}}\,(z^2-1)
-\frac{z^2}2\, \frac{\sqrt{1-\rho ^2}\,\psi''(\rho)}{\psi'(\rho)}. 
\end{equation*}
In particular, 
\begin{equation*}
\De_R(z)=\frac\rho{3\sqrt{1-\rho ^2}}\,(z^2-1), 
\end{equation*} 
\begin{equation*}
\De_{R_F}(z)=\De_{\psi_\infty(R)}(z)=-\frac\rho{3\sqrt{1-\rho ^2}}\,(2z^2+1), 
\end{equation*}
where $\psi_\infty$ is as defined in \eqref{eq:psi_infty}, and 
\begin{equation*}
\De_{\Psi_{4;\al}(R)}(z)=\De_{\psi_{4;z_\al}(R)}(z)
=\frac\rho{3\sqrt{1-\rho ^2}}\,\Big(\frac{z^2}{z_\al^2}-1\Big). 
\end{equation*} 

So, 
$|\De_{\Psi_{4;\al}(R)}(z_\be)|<|\De_{R_F}(z_\be)|$ for $\rho\ne0$ whenever $z_\al\ge1/\sqrt2$ or, equivalently, $\al\in(0,1-\Phi(1/\sqrt2))$, with $1-\Phi(1/\sqrt2)=0.2397\ldots$. 
Therefore, 
the transform $\Psi_{4;\al}$, which is asymptotically optimal for the given significance level $\al\in(0,1)$, will still be asymptotically better than Fisher's transform $R_F$ for any significance level $\be\in(0,0.5)$ provided that $\al\in(0,0.2397)$. 

As for the comparison of the asymptotically optimal transform $\Psi_\al(R)$ with $R$ itself, we can similarly see that, for instance, if $\al=0.05$, then $\Psi_\al(R)$ will be asymptotically better than $R$, not just for the significance level $\al=0.05$, but for any significance level $\be\in(0,0.11344)$; if $\al=0.01$, then $\Psi_\al(R)$ will be asymptotically better than $R$ for any significance level $\be\in(0,0.096927)$. 

We see that the advantage of the asymptotically optimal transform of $R$ over $R$ itself is substantially less in the SquareV model than in the BVN model. 
Vice versa, the advantage of the asymptotically optimal transform of $R$ over the Fisher transform $R_F$ of $R$ is much greater in the SquareV model than in the BVN model.

\subsection{Comparison of asymptotically optimal transform \texorpdfstring{$\Psi_{4;\al}(R)$}{} in the SquareV model with \texorpdfstring{$R$}{} and \texorpdfstring{$R_F$}{} for finite sample sizes} 
\label{finite-4}

Computer simulations 
mainly similar to those for finite sample sizes in the BVN model were carried out for the SquareV
model as well, the main difference with the BVN case being the use of the Mathematica command \texttt{RandomChoice[]} instead of 
\texttt{RandomVariate[BinormalDistribution[], \{\}]}. However, instead of about 3.6 hours for the BVN model, the simulations of the same number, 
$12\times10^6\times2\times4\times(10+10^2+10^3+10^4)\approx10^{12}$, of realizations of the random pair $(Y,Z)$ in the SquareV model took about 17 hours. It appears that the command \texttt{RandomChoice[]}, used to simulate distributions on finite sets, is not as highly optimized as  \texttt{RandomVariate[BinormalDistribution[], \{\}]}. 

After these simulations for the SquareV model had been completed, it was realized that the SquareV model can obtained by a transformation and a subsequent re-parametrization of the BVN model, as follows: if a random pair $(U,V)$ has the BVN distribution with zero means, unit variances, and an arbitrary correlation coefficient $\th\in(-1,1)$, then the random pair 
\begin{equation}\label{eq:BVN to 4pts}
(Y,Z):=(\sign U,\sign V) 	
\end{equation}
will have the SquareV distribution described by formula \eqref{eq:SM1}, provided that the parameters $\th$ and $\rho$ are related by the formula 
\begin{equation*}
(-1,1)\ni\cos[\tfrac\pi2\,(1-\rho)]=\th\longleftrightarrow\rho=1-\tfrac2\pi\,\arccos\th\in(-1,1). 	
\end{equation*}
This can be shown, for instance, by a standard change of variables in the relevant integrals to reduce the BVN density with correlation coefficient $\th$ to the standard BVN density, and then integrating in polar coordinates.  

So, the simulations for the SquareV model were mainly reduced to those for the BVN model, with the reduction of the execution time from about 17 hours to about 4 hours. Thereby, altogether $2\times12\times10^6\times2\times4\times(10+10^2+10^3+10^4)\approx2\times10^{12}$ realizations of the random pair $(Y,Z)$ in the SquareV model were obtained. 
Thus, one was able to reduce the corresponding estimates $\tilde s_{\cdot}(\al,\rho,n)$ of the standard error of the  estimates $\hat\vp_{\cdot}(\al,\rho,n)$ of the relative error of the approximation of $\al$ using $R$, $R_F$, and $\Psi_{4;\al}(R)$ by a factor of about $\sqrt2\approx1.4$.  

\begin{figure}[h]
	\centering
		\includegraphics[width=1.00\textwidth]{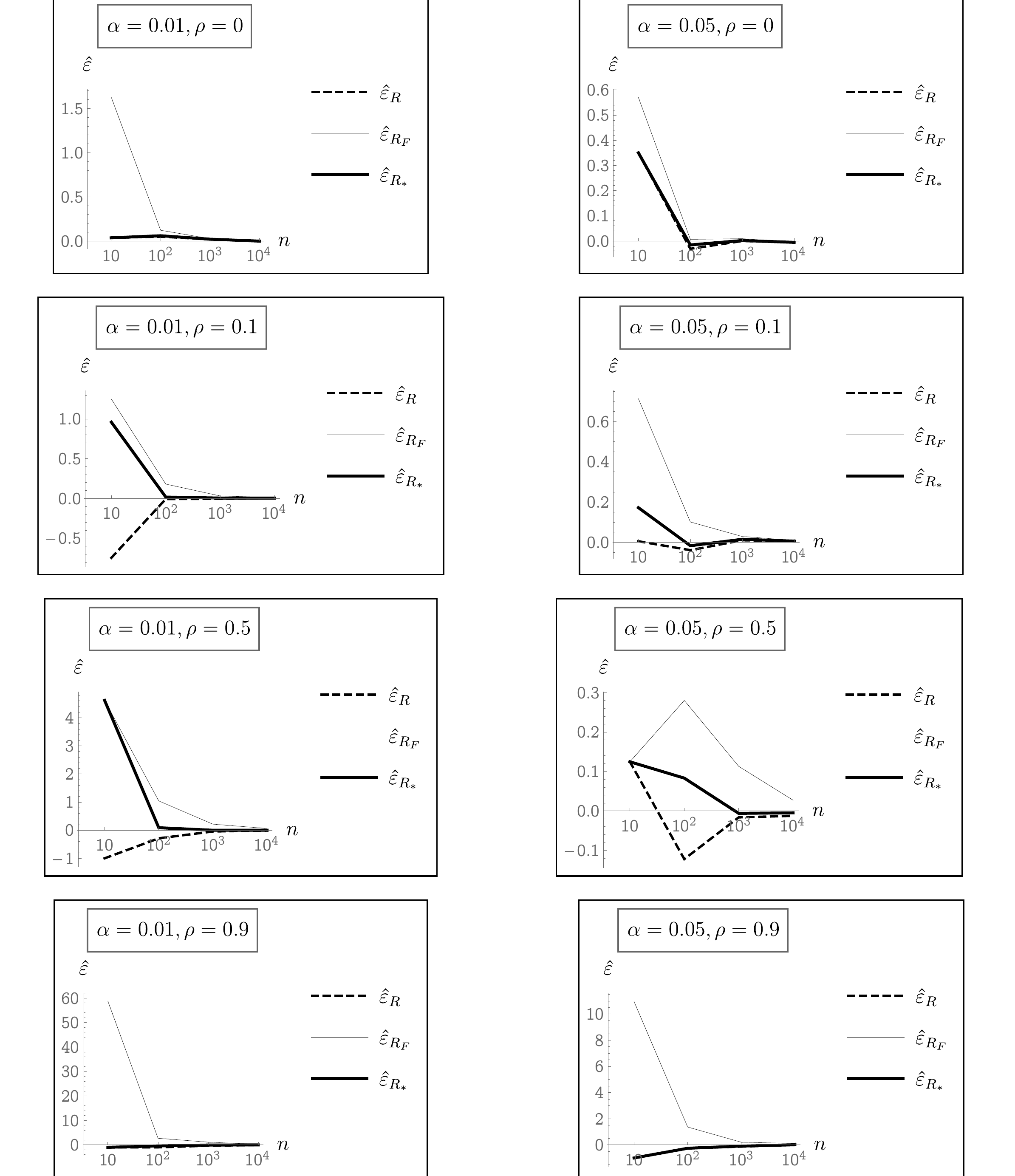}
	\caption{	Means 
	of the 
estimates 
of the relative errors of the approximation of the significance level $\al$ for Pearson's $R$ (thin), Fisher's transform $R_F$ of $R$ (dashed), and the asymptotically optimal transform $R_*:=\Psi_\al(R)$ (thick) in the SquareV model}
	\label{fig:pic1_4points}
\end{figure}

\begin{figure}[h]
	\centering
		\includegraphics[width=1.00\textwidth]{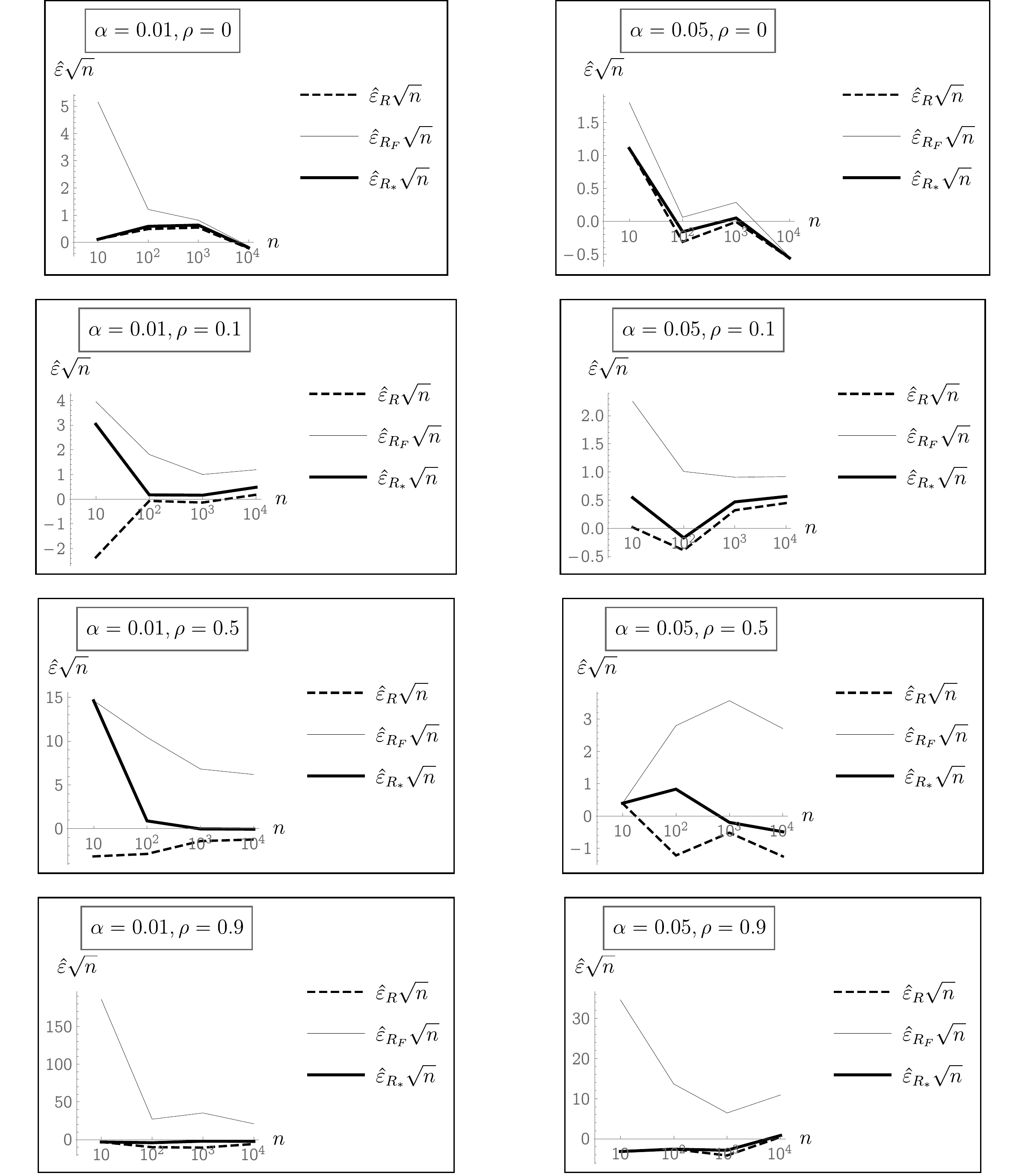}
	\caption{Means as in Fig.\ \ref{fig:pic1_4points} multiplied by $\sqrt n$}
	\label{fig:pic2_4points}
\end{figure}

Except for the differences described above, the format of Table~\ref{tab:tab4points} for the SquareV model is similar to that of Table~\ref{tab:tabBVN} for the BVN model. Correspondingly, Figs.~\ref{fig:pic1_4points} and \ref{fig:pic2_4points} for the SquareV model are similar to  Figs.~\ref{fig:pic1BVN} and \ref{fig:pic2BVN} for the BVN model.

\begin{figure}[h]
	\centering
		\includegraphics[width=1.00\textwidth]{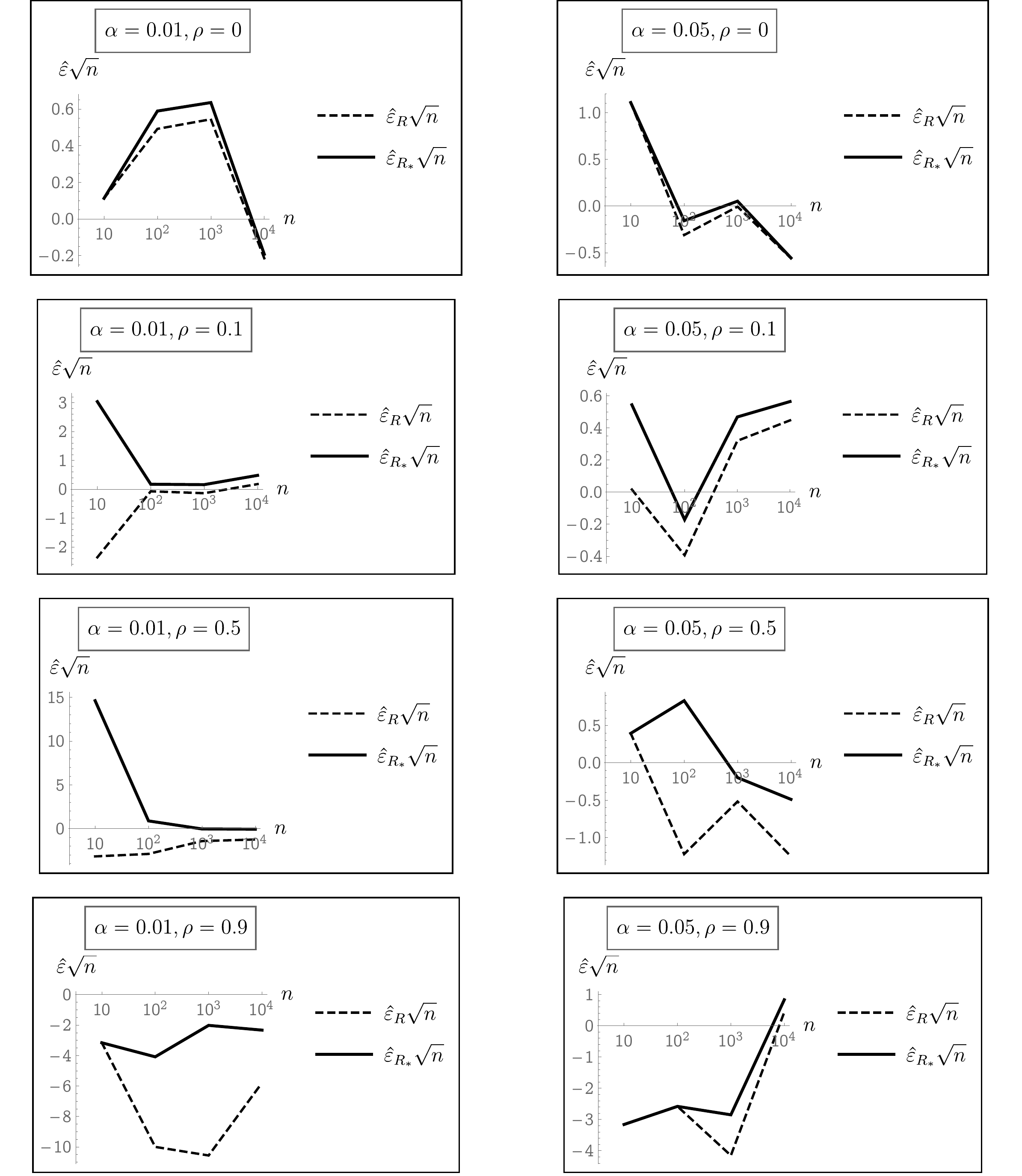}
	\caption{4points sqrt w/out Fisher's}
	\label{fig:pic2a_4points}
\end{figure}

Looking at Fig.~\ref{fig:pic2_4points}, we see that generally Fisher's transform $R_F$ performs much worse than $R$ and $\Psi_{4;\al}(R)$, at least when $\rho\in\{0.1,0.5,0.9\}$ and $n\in\{10^2,10^3,10^4\}$. Therefore, to make the comparisons between $R$ and $\Psi_{4;\al}(R)$ clearer, in Fig.~\ref{fig:pic2a_4points} the graphs corresponding to $R_F$ were removed. 

Somewhat similarly to the case of the BVN model, here we can see from Fig.~\ref{fig:pic2a_4points} that, 
if $n=10$ or if $\rho\in\{0,0.1\}$, the asymptotically optimal transform $\Psi_{4;\al}(R)$ appears to have comparatively little or no advantage over $R$ itsel;  
however, for $\rho\in\{0.5,0.9\}$ and $n\in\{10^2,10^3,10^4\}$, the transform $\Psi_{4;\al}(R)$ appears to be mostly better or much better than $R$ itself. 

\section{Conclusion}\label{concl} 
The main result of this paper is Theorem~\ref{th:general}, which shows for any correlation-parametrized model of dependence and for any given significance level $\al\in(0,1)$, there is an asymptotically optimal transform of Pearson's correlation statistic $R$, for which the generally leading error term for the normal approximation vanishes 
for all values $\rho\in(-1,1)$ of the correlation coefficient. 

It is also shown that in the BVN model Pearson's $R$ turns out to be asymptotically optimal for a rather unusual significance level $\al\approx0.240$, whereas Fisher's transform $R_F$ of $R$ is asymptotically optimal for the limit significance level $\al=0$. 
In the other specific model of dependence considered in this paper -- the SquareV model, Pearson's $R$ is asymptotically optimal for a still rather high significance level $\al\approx0.159$, whereas Fisher's transform $R_F$ of $R$ is not asymptotically optimal for any $\al\in[0,1]$. 

Moreover, we saw that in both the BVN model and the SquareV model, 
the transform $\Psi_\al(R)$, asymptotically optimal for a given value of $\al$, is in fact asymptotically better than $R$ and $R_F$ in wide ranges of values of the significance level, including $\al$ itself. 

Recall that Fisher's transform $R_F$ of $R$ was designed for the BVN case, with the purpose of making the asymptotic variance constant with respect to the correlation coefficient $\rho$. That $R_F$ usually turns out to be asymptotically closer to normality than $R$ in the BVN  model might now be explained by the observation that the significance level $\al=0$ (for which $R_F$ is asymptotically optimal in the BVN case) is closer to such usual in statistical practice values of the significance level as $0.05$ than to the significance level $\al_R\approx0.240$ (for which $R$ is asymptotically optimal in the BVN case). 

Extensive computer simulations for the BVN and SquareV models of dependence presented in this paper suggest that, for sample sizes $n\ge100$ and significance levels $\al\in\{0.01,0.05\}$, the mentioned asymptotically optimal transform of $R$ generally outperforms both Pearson's $R$ and Fisher's transform $R_F$ of $R$, the latter appearing generally much inferior to both $R$ and the asymptotically optimal transform of $R$ in the SquareV model.  

%
%


\begin{table}[h]
	\centering
{\small	
		\begin{tabular}{llr|lll}
\rule[-5pt]{0pt}{0pt} \multirow{2}{4ex}{$\al$}
& \multirow{2}{2ex}{$\rho$}
& \multirow{2}{5ex}{$n$}
& $\hat\vp_R(\al,\rho,n)$ & $\hat\vp_{R_F}(\al,\rho,n)$ & $\hat\vp_{\Psi_\al(R)}(\al,\rho,n)$ \\ 
& & & $\pm\tilde s_R(\al,\rho,n)$ & $\pm\tilde s_{R_F}(\al,\rho,n)$ & $\pm\tilde s_{\Psi_\al(R)}(\al,\rho,n)$ \\ 
\hline
\rule[12pt]{0pt}{0pt}0.01 & 0 & 10 & $-0.236\pm0.00253$ & $1.63\pm0.00463$ & $1.45\pm
   0.00439$ \\ 
 0.01 & 0 & 100 & $-0.00460\pm0.00344$ & $0.113\pm0.00359$ & $0.102\pm
   0.00353$ \\ 
 0.01 & 0 & 1000 & $-0.00306\pm0.00219$ & $0.00802\pm0.00204$ &
   $0.00689\pm0.00208$ \\ 
 0.01 & 0 & 10000 & $-0.00481\pm0.00262$ & $-0.00373\pm0.00263$ &
   $-0.00385\pm0.00264$ \\ 
 0.01 & 0.1 & 10 & $-0.683\pm0.00131$ & $1.72\pm0.00496$ & $1.47\pm
   0.00421$ \\ 
 0.01 & 0.1 & 100 & $-0.139\pm0.00355$ & $0.123\pm0.00375$ & $0.0971\pm
   0.00365$ \\ 
 0.01 & 0.1 & 1000 & $-0.0432\pm0.00263$ & $0.0133\pm0.00279$ &
   $0.00813\pm0.00285$ \\ 
 0.01 & 0.1 & 10000 & $-0.0128\pm0.00256$ & $0.00265\pm0.00245$ &
   $0.00114\pm0.00248$ \\ 
 0.01 & 0.5 & 10 & $-1.00\pm0$ & $2.08\pm0.00419$ & $1.53\pm0.00396$ \\ 
 0.01 & 0.5 & 100 & $-0.562\pm0.00161$ & $0.180\pm0.00237$ & $0.0980\pm
   0.00245$ \\ 
 0.01 & 0.5 & 1000 & $-0.200\pm0.00274$ & $0.0298\pm0.00295$ & $0.00697\pm
   0.00278$ \\ 
 0.01 & 0.5 & 10000 & $-0.0638\pm0.00329$ & $0.00837\pm0.00391$ &
   $0.00134\pm0.00385$ \\ 
 0.01 & 0.9 & 10 & $-1.00\pm0$ & $2.41\pm0.00642$ & $1.56\pm0.00694$ \\ 
 0.01 & 0.9 & 100 & $-0.833\pm0.00109$ & $0.236\pm0.00324$ & $0.0965\pm
   0.00280$ \\ 
 0.01 & 0.9 & 1000 & $-0.342\pm0.00270$ & $0.0436\pm0.00273$ & $0.00382\pm
   0.00295$ \\ 
 0.01 & 0.9 & 10000 & $-0.119\pm0.00279$ & $0.00877\pm0.00290$ &
   $-0.00352\pm0.00281$ \\ 
 0.05 & 0 & 10 & $0.233\pm0.00115$ & $0.626\pm0.00199$ & $0.559\pm
   0.00182$ \\ 
 0.05 & 0 & 100 & $0.0217\pm0.00121$ & $0.0530\pm0.00131$ & $0.0473\pm
   0.00133$ \\ 
 0.05 & 0 & 1000 & $0.00309\pm0.00124$ & $0.00609\pm0.00120$ & $0.00554\pm
   0.00121$ \\ 
 0.05 & 0 & 10000 & $0.00131\pm0.00121$ & $0.00167\pm0.00122$ &
   $0.00160\pm0.00122$ \\ 
 0.05 & 0.1 & 10 & $0.0395\pm0.00140$ & $0.675\pm0.00161$ & $0.569\pm
   0.00166$ \\ 
 0.05 & 0.1 & 100 & $-0.0292\pm0.00129$ & $0.0597\pm0.00143$ & $0.0432\pm
  0.00145$ \\ 
 0.05 & 0.1 & 1000 & $-0.0131\pm0.00150$ & $0.00772\pm0.00143$ &
   $0.00383\pm0.00143$ \\ 
 0.05 & 0.1 & 10000 & $-0.00505\pm0.000957$ & $0.000768\pm0.00106$ &
   $-0.000367\pm0.00104$ \\ 
 0.05 & 0.5 & 10 & $-0.723\pm0.000849$ & $0.861\pm0.00128$ & $0.587\pm
   0.00143$ \\ 
 0.05 & 0.5 & 100 & $-0.218\pm0.000924$ & $0.103\pm0.00123$ & $0.0435\pm
   0.00111$ \\ 
 0.05 & 0.5 & 1000 & $-0.0703\pm0.000603$ & $0.0219\pm0.000734$ &
   $0.00470\pm0.000712$ \\ 
 0.05 & 0.5 & 10000 & $-0.0224\pm0.00119$ & $0.00578\pm0.00124$ &
   $0.000620\pm0.00120$ \\ 
 0.05 & 0.9 & 10 & $-1.00\pm5.75\times10^{-6}$ & $1.04\pm0.00193$ &
   $0.588\pm0.00160$ \\ 
 0.05 & 0.9 & 100 & $-0.403\pm0.000995$ & $0.145\pm0.00168$ & $0.0419\pm
   0.00156$ \\ 
 0.05 & 0.9 & 1000 & $-0.128\pm0.000776$ & $0.0347\pm0.000863$ &
   $0.00463\pm0.000964$ \\ 
 0.05 & 0.9 & 10000 & $-0.0425\pm0.00124$ & $0.00825\pm0.00118$ &
   $-0.00118\pm0.00122$ \\ 
		\end{tabular}
		}
		\medskip
	\caption{\rule[12pt]{0pt}{0pt}Means and standard deviations of the 
estimates \eqref{eq:hateps} of the relative errors of the approximation of the significance level $\al$ for Pearson's $R$, Fisher's transform $R_F$ of $R$, and the asymptotically optimal transform $\Psi_\al(R)$ in the BVN model}
	\label{tab:tabBVN}
\end{table}

\begin{table}[h]
	\centering
{\small	
		\begin{tabular}{llr|lll}
\rule[-5pt]{0pt}{0pt} \multirow{2}{4ex}{$\al$}
& \multirow{2}{2ex}{$\rho$}
& \multirow{2}{5ex}{$n$}
& $\hat\vp_R(\al,\rho,n)$ & $\hat\vp_{R_F}(\al,\rho,n)$ & $\hat\vp_{\Psi_{4;\al}(R)}(\al,\rho,n)$ \\ 
& & & $\pm\tilde s_R(\al,\rho,n)$ & $\pm\tilde s_{R_F}(\al,\rho,n)$ & $\pm\tilde s_{\Psi_{4;\al}(R)}(\al,\rho,n)$ \\ 
\hline
\rule[12pt]{0pt}{0pt}0.01 & 0 & 10 & $0.0356\pm 0.00173$ & $1.63\pm 0.00329$ & $0.0356\pm
   0.00173$ \\
 0.01 & 0 & 100 & $0.0492\pm 0.00215$ & $0.121\pm 0.00226$ & $0.0589\pm
   0.00220$ \\
 0.01 & 0 & 1000 & $0.0172\pm 0.00238$ & $0.0257\pm 0.00232$ & $0.0201\pm
   0.00236$ \\
 0.01 & 0 & 10000 & $-0.00215\pm 0.00190$ & $-0.00133\pm 0.00188$ &
   $-0.00193\pm 0.00190$ \\
 0.01 & 0.1 & 10 & $-0.746\pm 0.00101$ & $1.25\pm 0.00246$ & $0.959\pm
   0.00236$ \\
 0.01 & 0.1 & 100 & $-0.00727\pm 0.00172$ & $0.181\pm 0.00216$ & $0.0172\pm
   0.00185$ \\
 0.01 & 0.1 & 1000 & $-0.00439\pm 0.00170$ & $0.0315\pm 0.00178$ &
   $0.00500\pm 0.00173$ \\
 0.01 & 0.1 & 10000 & $0.00177\pm 0.00229$ & $0.0119\pm 0.00225$ &
   $0.00480\pm 0.00227$ \\
 0.01 & 0.5 & 10 & $-1.00\pm 0$ & $4.62\pm 0.00491$ & $4.62\pm 0.00491$ \\
 0.01 & 0.5 & 100 & $-0.288\pm 0.00197$ & $1.04\pm 0.00331$ & $0.0887\pm
   0.00237$ \\
 0.01 & 0.5 & 1000 & $-0.0447\pm 0.00190$ & $0.215\pm 0.00205$ &
   $-0.000958\pm 0.00200$ \\
 0.01 & 0.5 & 10000 & $-0.0123\pm 0.00177$ & $0.0618\pm 0.00182$ &
   $-0.000646\pm 0.00182$ \\
 0.01 & 0.9 & 10 & $-1.00\pm 0$ & $58.8\pm 0.00975$ & $-1.00\pm 0$ \\
 0.01 & 0.9 & 100 & $-1.00\pm 0$ & $2.71\pm 0.00477$ & $-0.409\pm 0.00149$
   \\
 0.01 & 0.9 & 1000 & $-0.334\pm 0.00138$ & $1.12\pm 0.00256$ & $-0.0639\pm
   0.00172$ \\
 0.01 & 0.9 & 10000 & $-0.0579\pm 0.00214$ & $0.209\pm 0.00220$ &
   $-0.0233\pm 0.00225$ \\
 0.05 & 0 & 10 & $0.350\pm 0.000936$ & $0.570\pm 0.00103$ & $0.350\pm
   0.000936$ \\
 0.05 & 0 & 100 & $-0.0311\pm 0.000923$ & $0.00629\pm 0.000933$ &
   $-0.0157\pm 0.000939$ \\
 0.05 & 0 & 1000 & $-0.000207\pm 0.000941$ & $0.00914\pm 0.000940$ &
   $0.00163\pm 0.000936$ \\
 0.05 & 0 & 10000 & $-0.00559\pm 0.000614$ & $-0.00538\pm 0.000611$ &
   $-0.00556\pm 0.000612$ \\
 0.05 & 0.1 & 10 & $0.00496\pm 0.000823$ & $0.714\pm 0.00107$ & $0.172\pm
   0.000941$ \\
 0.05 & 0.1 & 100 & $-0.0393\pm 0.000888$ & $0.101\pm 0.00106$ & $-0.0172\pm
   0.000912$ \\
 0.05 & 0.1 & 1000 & $0.0101\pm 0.000757$ & $0.0287\pm 0.000741$ &
   $0.0148\pm 0.000740$ \\
 0.05 & 0.1 & 10000 & $0.00447\pm 0.000839$ & $0.00916\pm 0.000829$ &
   $0.00564\pm 0.000832$ \\
 0.05 & 0.5 & 10 & $0.125\pm 0.00110$ & $0.125\pm 0.00110$ & $0.125\pm
   0.00110$ \\
 0.05 & 0.5 & 100 & $-0.122\pm 0.000917$ & $0.280\pm 0.000882$ & $0.0831\pm
   0.000959$ \\
 0.05 & 0.5 & 1000 & $-0.0164\pm 0.00101$ & $0.113\pm 0.000981$ &
   $-0.00627\pm 0.00102$ \\
 0.05 & 0.5 & 10000 & $-0.0125\pm 0.000649$ & $0.0271\pm 0.000660$ &
   $-0.00489\pm 0.000682$ \\
 0.05 & 0.9 & 10 & $-1.00\pm 0$ & $11.0\pm 0.00267$ & $-1.00\pm 0$ \\
 0.05 & 0.9 & 100 & $-0.258\pm 0.000722$ & $1.37\pm 0.00112$ & $-0.258\pm
   0.000722$ \\
 0.05 & 0.9 & 1000 & $-0.131\pm 0.000720$ & $0.204\pm 0.000717$ &
   $-0.0902\pm 0.000712$ \\
 0.05 & 0.9 & 10000 & $0.00435\pm 0.000870$ & $0.109\pm 0.000945$ &
   $0.00831\pm 0.000871$ \\
		\end{tabular}
		}
		\medskip
	\caption{\rule[12pt]{0pt}{0pt}4points}
	\label{tab:tab4points}
\end{table}

\bibliographystyle{imsart-number}


\bibliography{P:/pCloudSync/mtu_pCloud_02-02-17/bib_files/citations10.13.18a}

\end{document}